\documentclass[]{interact}

\usepackage{epstopdf}
\usepackage[caption=false]{subfig}

\usepackage[numbers,sort&compress]{natbib}
\bibpunct[, ]{[}{]}{,}{n}{,}{,}
\makeatletter
\def\NAT@def@citea{\def\@citea{\NAT@separator}}
\makeatother

\theoremstyle{plain}

\theoremstyle{definition}

\theoremstyle{remark}

\begin{document}


\title{Multiagent Control of Airplane Wing Stability with ``Feathers'' under the Flexural Torsional Flutter}

\author{
\name{Dmitry.~S. Shalymov\textsuperscript{a,b}, Oleg N. Granichin\textsuperscript{a,b},
Yury Ivanskiy\textsuperscript{a,b} and Zeev Volkovich\textsuperscript{c}
\thanks{CONTACT Y. Ivanskiy. Email: y.ivansky@spbu.ru} 
}
\affil{
\textsuperscript{a}St. Petersburg State University, Universitetskaya nab. 7-9, Saint-Petersburg, 199034, Russia; \\
\textsuperscript{b}IPME RAS, V.O., Bolshoj pr., 61, St. Petersburg, 199178, Russia;\\
\textsuperscript{c}Software Engineering Department, Ort Braude College, Rehov Snunit 51, POB 78, Karmiel 2161002, Israel
}
}

\maketitle

\begin{abstract}
This paper proposes a novel method for the prevention of the unbounded oscillation of an aircraft wings under the flexural torsion flutter. The paper introducing the novel multiagent method for control of an aircraft wing, assuming that the wing surface consists of controlled ``feathers'' (agents). Theoretical evaluation of the approach demonstrates its high ability to prevent flexural-torsional vibrations of an aircraft. Since the model expands the possibilities for damping the wing oscillations, which potentially allows an increase in aircraft speed without misgiving of flutter. The study exhibits that timing is the main limitation of dampening vibrations. A new method for controlling an aircraft wing is suggested to increase the maximal flight speed of an aircraft without flutter occurrence via a novel model of the bending-torsional vibrations of an airplane wing with controlled feathers on its surface and new control laws based on the Speed-Gradient methodology.  Provided experiments demonstrate the theoretical and advantages of the multiagent approach to ``feathers'' movement control.
\end{abstract}


\section{INTRODUCTION}

The paper investigates a critical known problem of airplane wing control under the flexural torsion flutter. The stability in this task is closely connected to the ability to maintain under oscillations wings integrity because, during a flight, their components can experience significant strain. After attaining a certain flight speed, say $V_{flat}$, oscillations of the wing rapidly and catastrophically increase until the wing breaks if they are not stopped. It is so named the flutter phenomenon. For example, even in horizontal flight at a constant speed on a heavy transport aircraft, the deflection of the wing end can achieve several meters~\cite{A1} so that the corresponding deformations affect the magnitude and distribution of the aerodynamic load leading potentially to structural instability, both static (e.g. wing divergence) and dynamic (e.g. flutter).

Therefore, safeguarding the required aerodynamic characteristics and wing stability in various aircraft flight phases are needed. It should be noted that the elements of the wing mechanization are widespread in modern aircraft construction (pre-flaps, ailerons, flaps, etc.) serving precisely this purpose, in particular, in the most critical take-off and landing modes. However, there is no current common way to counteract wing flutter effectively. Maybe except to impose a bound on the aircraft maximum flight speed $V_{max}$ to be less than $V_{flat}$ 
Inventors strive to enlarge flight speed $V_{flat}$ as much as possible in order to increase the maximum allowable aircraft speed $V_{max}$.

\begin{figure}[thpb]
	\centering
	\includegraphics[scale=0.4]{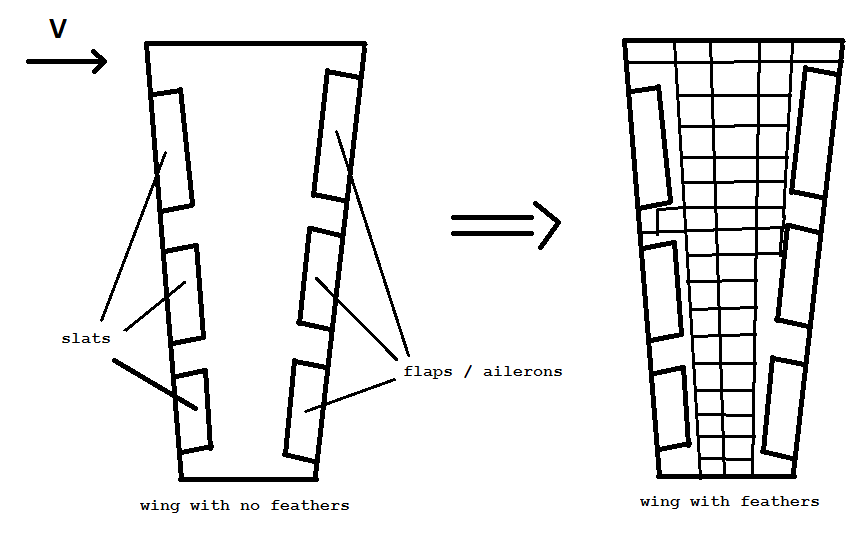}
	\caption{Wing with feathers}
	\label{fig:1p}
\end{figure}

A novel approach proposed in this paper suggests to cover an aircraft wing with small-sized movable elements (``feathers'') both above and below, capable of changing their orientation consistently with the airflow. In the neutral position, when feathers are not raised but lie on the surface, they do not affect the calculated wing profile (see Fig.~\ref{fig:1p}, the left panel). In this case the wing dynamic does not deteriorate.

The methodology novelty lies in the use of ``feathers'' on the aircraft's surface. In this way, a wing becomes a multiagent system that can be controlled via the control theory methods. 
At the heart of multiagent systems is a decentralized approach to solving problems, in
which dynamically updated information in a distributed network of
intelligent agents is processed not in some center, but directly at the
agents on the basis of their local observations and locally available
information from neighbors. At the same time, both resource and time
costs for communication in the network are significantly reduced, as well as
the time for processing and decision-making in the center of the entire
system (if it does exist).
Such an approach is new and has not been considered before in the flutter problem area. Because of the relatively intensive wings oscillations during flutter, their essential performance indicator is embodied by a period needed to damp and maintain vibrations safely. 
Application of multiagent approach allows an occurrence of emergent intelligence (intellectual resonance, swarm intelligence) or an appearance of unexpected properties that a system possesses, but none of its individual elements has.
Each feather tries to minimize the deviation of the wing segment, to which it is attached, from its initial state. The action of each feather, in general, is not consolidated with other feathers actions, however a combination of all feathers impact results in a new property of the wing as the multiagent system of feathers to damp the vibrations. 

A crucial theoretical suggestion is the ability of a ``feather'' to swap the orientation instantly along with the airflow. The impact of a possible delay should be a subject of a further research.

The rest of the paper is structured in the following way: Section 2 provides an overview of related works. Section 3 is devoted to studying system dynamics for small flexural-torsional vibrations. Section 4 describes the control laws produced using the Speed-Gradient method. Section 5 presents a possible connection to the multiagent control methodology. Section 6 is devoted to the conclusion.

\section{RELATED WORK}

The condition frontier between stable and self-sustaining motions in a flight is named flutter speed and flutter boundary~\cite{WrightCopper}. The flutter can be divided into several groups according to the instability appearing with changes in conditions like dynamic pressure increasing. An explosive flutter occurs after exceeding the flutter speed $V_{flat}$. This process results in highly divergent oscillations and wing breaks within a fraction of a second.

The moderate (mild) flutter corresponds to when the system is stable but lowly significantly damped before achieving the flutter speed and can be identified well below the flutter speed by instabilities extrapolating.

An approach intended to steady an unstable flutter system is called Active Flutter 
Suppression (AFS). 
A broad overview of AFS research is presented in~\cite{G1}. 
The study of the ability to suppress flutter instability through actively controlled closed-loop action of control surfaces has a long history~\cite{c5}. 
Numerous researches in this field are carried out as early as the 1970s and 1980s~\cite{c9}. 
Various approaches to synthesizing AFS systems control law, including adaptive control methods and control with variable parameters, are also considered~\cite{c364,c366,c369,c370}. 
AFS is essential for an effective solution of aeroelastic instability problems and can lead to significant aircraft and airframes weight savings, (see, e.g.~\cite{c6}). 

The AFS approach based on elimination of delays in loading growth induced by unsteady aerodynamic stresses is investigated in~\cite{c371}. 
The concepts of ``active flexible wing'' or ``active flexible airframe'' are respected in~\cite{c57}; similarly, flight stability and controllability of rigid and flexible aircraft are considered in~\cite{c15, c23, c27}. The influence of the aeroelasticity on stability and controllability of flight using corrections of the derivatives of static aeroelastic stability is studied in~\cite{c31}; an active controls perspective is presented in~\cite{c134}.

The idea of using an active control system has been considered and discussed since the advent of flights with a human crew~\cite{c15,c23}. The adaptive control methodology appears attractive due to multiple plant characteristics and the ability to respond to damage scenarios, cf.~\cite{c417,c418,c420,c422,c423,c424}---article~\cite{c403} attempts to systematize the modern control theory laws developed in the field. A broad review of methods for synthesizing control laws and modeling of aeroelastic systems is presented in~\cite{c444}. With this connection, it is possible to recall several works considering different aspects of the mentioned thematics~\cite{c434, c437, c438, c439, c443},~\cite{c366,c369,c370}.

Multiagent systems have many applications in civilian, security, and military areas~\cite{cdc2013, ecc2015}. Centralized quantitative and qualitative modeling, analysis, constraint satisfaction, maintenance, and control seem to be too strict for these systems~\cite{math2021}. On the other hand, the distributed and incremental reasoning seems to be more scalable, robust and flexible. That is why an investigation of multiagent control systems is popular nowadays~\cite{F1, Gr3, Fradkov1}. The multiagent methodology can serve the general model of the interactions in a complex system~\cite{F7, Fradkov1}. An overview of publications considering the emergent intelligence and self-organization in groups of devices is provided in~\cite{math2021}.

Dynamical networks and multiagent protocol for the airplane are investigated in~\cite{Gran1,Gran2,A5}, where the multiagent control is used for leveling the perturbing forces in the turbulence.

\section{DYNAMICS EQUATIONS OF WING WITH FEATHERS}

\begin{figure}[thpb]
	\centering
	\includegraphics[scale=0.3]{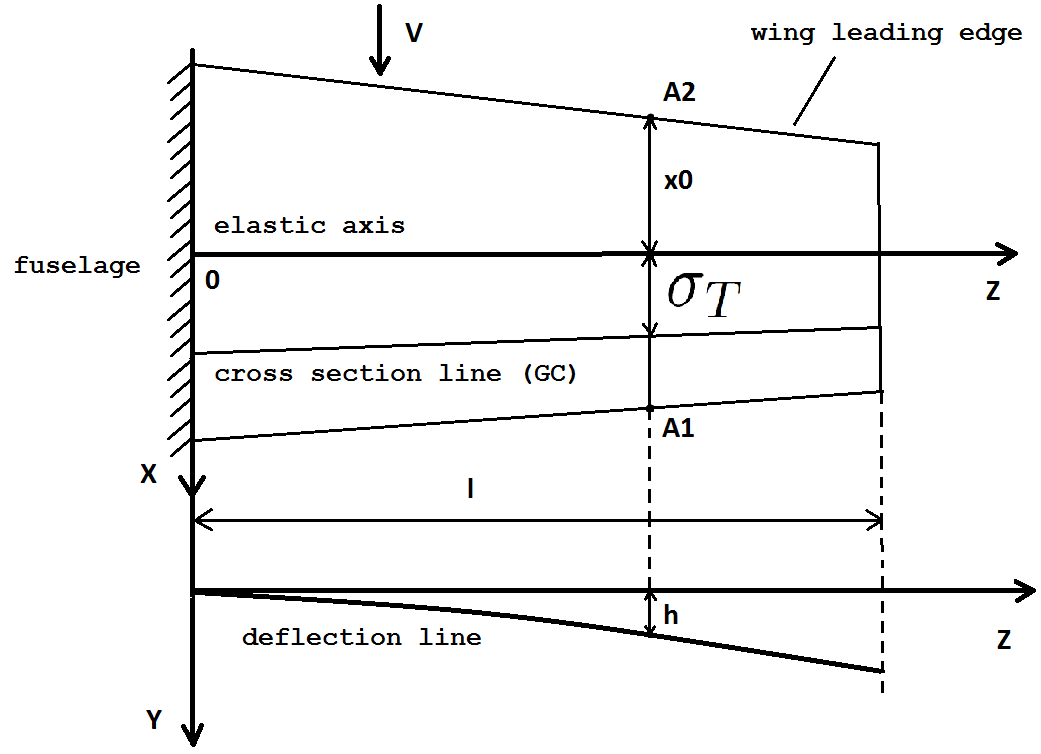}
	\caption{\centering
	Stiffness Centers and Gravity Centers lines on a wing, where 
	\parbox{0.65\textwidth}{
	$O$ is the stiffness center of the wing section on the fuselage; \\
	$X$-axis is directed on a free stream; \\
	$Z$-axis is directed along the elastic axis of the undeformed wing; \\
	$Y$-axis complements the coordinate system to the right; \\
	$x_{0}$ is the distance from the leading edge of the wing to the SC section; \\
	$\sigma_{T}$ is the distance between SC and GC; \\
	$h$ is the transverse deflection of the section $A_{1}A_{2}$; \\
	$l$ is the half-span of the wing.}
	}
	\label{fig:2p}
\end{figure}

This section describes the proposed approach named a ``wing with feathers''. Let us consider the phenomenon of flexural-torsional flutter as described in~\cite{A2,A3} in a steady horizontal flight at a constant speed. We regard the non-sweeping wing of the half-span $l$ and feathers in neutral position with a cantilevered beam with a static distributed load on bending and torsion.
The elastic axis of the wing passes through the Stiffness Centers (SC's) of the sections and does not coincide with the line of Gravity Centers (GC's) of the sections. We assume that the wing stiffness in the longitudinal and transverse directions of the wing plane to be very large neglecting of the vibrations in these directions. We also ignore the possible movements of the SC and GC along the sections during the flight. Figure~\ref{fig:2p} illustrates the conceivable location of the SC  and GC lines on a wing.

It must be noted that typically the GC is located behind the SC section.

\begin{figure}[thpb]
	\centering
	\includegraphics[scale=0.5]{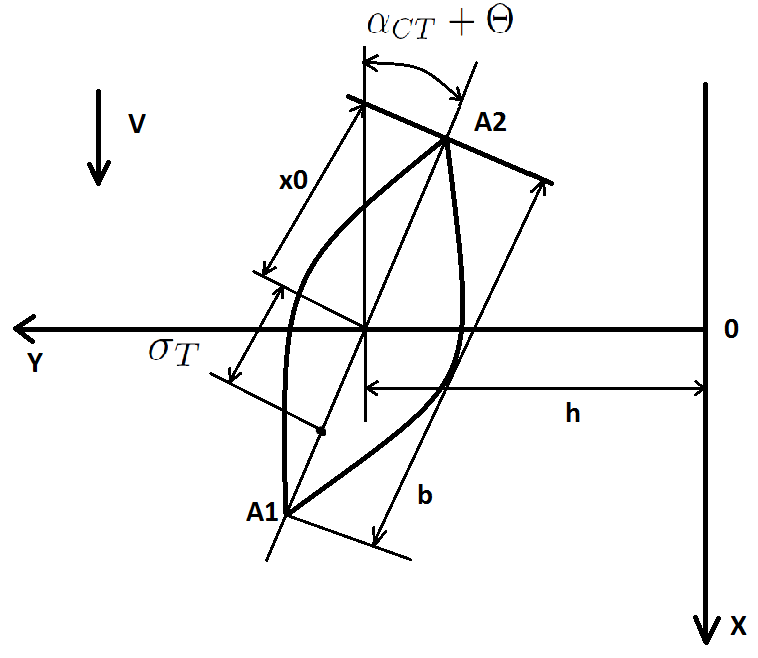}
	\caption{Wing cross section}
	\label{fig:3p}
\end{figure}


The corresponding equations of the elastic line of the beam are of the form, cf. \cite{c5}:

\begin{equation}\label{eq1}
\begin{cases}
\frac{\partial^{2}}{\partial z^{2}}\left( EJ\frac{\partial^{2}y}{\partial z^{2}}\right) = q^{0},\\
\frac{\partial}{\partial z}\left( GJ_{k}\frac{\partial\Theta}{\partial z}\right) = m^{0},
\end{cases}
\end{equation}
where as presented in Fig.~\ref{fig:3p}
\begin{itemize}
	\item $b$ is the wing section chord;
	\item $h$ is the transverse deflection of the section $A_{1}A_{2}$;
	\item $\alpha_{CT}$ is the angle of attack in the section of the undeformed wing;
	\item $y$ is the deflection of the stiffness axis in the current wing section;
	\item $\Theta$ is the wing twist angle, which is considered positive if it increases the angle of attack in the section;
	\item $EJ, GJ_{k}$ are the wing stiffness in bending and torsion, respectively;
	\item $q^{0}$ and $m^{0}$ are the linear force and moment relative to the stiffness axis acting on the wing.
\end{itemize}

The functions of \eqref{eq1} are time-independent since the wing is in a steady (stationary) state:
\begin{equation}\label{eq2}
y=y^{0}(z),~\Theta=\Theta^{0}(z).
\end{equation}

These solutions must satisfy the boundary conditions at the ends of the wing~\cite{fung}
\begin{equation}\label{eq3}
\begin{cases}
\left. y\right\rvert_{z=0}=0; ~~ \left. \frac{\partial y}{\partial z}\right\rvert_{z=0}=0~~ (\textnormal{tight fuselage fit});\\

\left. EJ\frac{\partial^{2}y}{\partial z^{2}}\right\rvert_{z=l}=0;~~ 
\left. \frac{\partial}{\partial z}\left(EJ\frac{\partial^{2}y}{\partial z^{2}}\right) \right\rvert_{z=l}=0 \\
~~ (\textnormal{moment and shear force at the free end}),\\

\left. \Theta\right\rvert_{z=0}=0; ~~ \left. GJ_{k}\frac{\partial \Theta}{\partial z}\right\rvert_{z=l}=0 ~~ \\
(\textnormal{angle of rotation of the terminated end and moment at the free end}).
\end{cases}
\end{equation}

Now, suppose that for some unknown reason like sudden aileron movement, air holes, a gust of wind, and so on, the wing deviated from its stationary position (see \eqref{eq2}). At the end of the specified impact, after the cessation of this factor exposure, the wing returns to the equilibrium state under
 the influence of the elastic forces. If the energy dissipation is insignificant, the aperiodic process does not appear. However, wing oscillations may arise. We assume that these fluctuations are initially ignored and do not affect the aircraft dynamics.

According to \cite{A2,A3}, small bending-torsional oscillations of a wing near its equilibrium position (see \eqref{eq2}, \eqref{eq3}) in the laminar flow are described by the following equations:
\begin{equation}\label{eq4}
\begin{cases}
\frac{\partial^{2}}{\partial z^{2}}\left( EJ\frac{\partial^{2}y_{1}}{\partial z^{2}}\right) + m\frac{\partial^{2}y_{1}}{\partial t^{2}} - m\sigma_{T}\frac{\partial^{2}\Theta_{1}}{\partial t^{2}}=q_{a}, \\

\frac{\partial}{\partial z}\left( GJ_{k}\frac{\partial\Theta_{1}}{\partial z}\right) + m\sigma_{T}\frac{\partial^{2}
y_{1}}{\partial t^{2}} - J_{m}\frac{\partial^{2}\Theta_{1}}{\partial t^{2}}=m_{a}, 
\end{cases}
\end{equation}

where 
\begin{itemize}
	\item $y_{1} $ and $\theta_{1}$ are the additional deflection and angle of twisting of the wing relative to the stationary state (see \eqref{eq2},\eqref{eq3})), due to fluctuations;
	\item $m$ is the linear mass of the wing;
	\item $J_{m}$ is the linear mass moment of inertia of the wing relative to its stiffness axis;
	\item $q_ {a}$ and $m_ {a}$ are the linear aerodynamic force of the wing and the linear moment of the aerodynamic force relative to the stiffness axis, due to wing vibrations.
\end{itemize}

Solutions of the system \eqref{eq4} must satisfy boundary conditions similar to \eqref{eq3}. We represent the right-hand sides of \eqref{eq4} in the form:
\begin{equation*}
q_{a}=\Delta q_{a} + q_{u}, ~~ m_{a} = \Delta m_{a} + m_{u},
\end{equation*}

where 
\begin{itemize}
	\item $\Delta q_ {a}$ and $\Delta m_ {a}$ are the linear aerodynamic force and moment relative to the stiffness axis respectively, arising due to wing oscillations in the neutral position of the feathers;
	\item $q_{u}$ and $m_{u}$ are the linear aerodynamic force and moment created by changing the orientation of the feathers.
\end{itemize}

Following the flutter equations (\cite{A2}, p. 176, equation (35)) and taking into account linear aerodynamic forces and moment $q_{u}$ and $m_{u}$, we rewrite \eqref{eq4} as

\begin{equation}\label{sist6}
\begin{cases}
\frac{\partial^{2}}{\partial z^{2}}(EJ\frac{\partial^{2}y_{1}}{\partial z^{2}})
+m\frac{\partial^{2}y_{1}}{\partial t^{2}}
-m\sigma_{T}\frac{\partial^{2}\Theta_{1}}{\partial t^{2}} \\
-C_{y}^{\alpha}\left[\Theta_{1}+(\frac{3}{4}b-x_{0})\frac{1}{V}\frac{\partial\Theta_{1}}{\partial t} -\frac{1}{V}\frac{\partial y_{1}}{\partial t}\right]\rho bV^{2}=q_{u}\\

\frac{\partial}{\partial z}(GJ_{k}\frac{\partial\Theta_{1}}{\partial z}) + m\sigma_{T}\frac{\partial^{2}y_{1}}{\partial t^{2}} - J_{m}\frac{\partial^{2}\Theta_{1}}{\partial t^{2}} 
- \frac{\pi}{16}\frac{b^2}{V}\frac{\partial\Theta_{1}}{\partial t} \rho bV^{2} \\
+\left\{
+C_{y}^{\alpha}(x_{0}-\frac{b}{4})\left[\Theta_{1}+(\frac{3}{4}b-x_{0})\frac{1}{V}\frac{\partial\Theta_{1}}{\partial t} -\frac{1}{V}\frac{\partial y_{1}}{\partial t}\right]
\right\}\rho bV^{2} = m_{u}, \\

y_{1} = \frac{\partial y_{1}}{\partial z}=\Theta_{1} = 0,~~z=0,\\

\frac{\partial^{2}y_{1}}{\partial z^{2}} = \frac{\partial^{3}y_{1}}{\partial z^{3}} = \frac{\partial\Theta_{1}}{\partial z}=0,~~z=l,

\end{cases}
\end{equation}

where 
\begin{itemize}
	\item $C_{y}^{\alpha}=\frac{\partial C_{y}}{\partial \alpha}$; $C_{y}$ is the wing lift coefficient;
	\item $C_{y}^{\alpha}$ consider constant along the span;
	\item $C_{y}=C_{y}^{\alpha}(\alpha-\alpha_{0});$
	\item $\alpha=\alpha_{CT}+\Theta^{0}+\Theta_{1}$ is the instant value of the angle of attack when the wing moves;
	\item $\alpha_{0}$ is the value of the angle of attack at which $C_{y}=0$;
	\item $\rho$ is the air density.
\end{itemize}

As it can be seen from \eqref{sist6}, the bending and torsional vibrations of the wing are interdependent. It is one of the necessary conditions for flutter occurrence. 
It is also known \cite{c5} with increasing speed $V$, the bending and torsional oscillations approach each other, and for $V = V_{flat}$ a wing coalesce. Moreover, there is a phase shift between these oscillations, a necessary condition for the occurrence of flutter \cite{A2,A3}.

It is important to remark that in this case the wing amplitude oscillations wavers around a small constant value. So the oscillations themselves are no longer self-damped, as it was the case where $V < V_{flat}$. The crucial problem arises when $V > V_{flat}$, i.e. when the slightest deformations overgrow catastrophically. 

In terms of \eqref{sist6}, it is necessary to form $q_u$ and $m_u$ to prevent flutter so that a wing oscillation energy is bounded. The bound value is reliable from a controllability standpoint, taking into account the stability and the integrity of the aircraft structure. So, if $V > V_{flat}$, for the total energy of the cantilevered beam, \cite{A2}, the condition is:

\begin{multline}\label{eq7}
E = E_{kinet} + E_{poten} = 
\frac{1}{2}\int_{0}^{l}{m\left(\frac{\partial y_{1}}{\partial t} \right)^{2}dz}
+ \frac{1}{2}\int_{0}^{l}{J_{m}\left(\frac{\partial \Theta_{1}}{\partial t} \right)^{2}dz} 
-\int_{0}^{l}{m\sigma_{T}\frac{\partial y_{1}}{\partial t}\frac{\partial \Theta_{1}}{\partial t}dz} + \\ 
\frac{1}{2}\int_{0}^{l}{EJ\left(\frac{\partial^{2}y_{1}}{\partial z^{2}} \right)^{2}dz} 
+\frac{1}{2}\int_{0}^{l}{GJ_{k}\left(\frac{\partial \Theta_{1}}{\partial z} \right)^{2}dz} \leq E_{*}.
\end{multline}

A stricter requirement is to get the system \eqref{sist6} solution in a given small neighborhood of the solution of \eqref{eq2}:
\begin{equation}\label{eq8}
\left\|\bar{x}\right\| \leq \epsilon,
\end{equation}
where 
\begin{equation}\notag
\bar{x} = 
\left( y_{1}, \frac{\partial y_{1}}{\partial t}, \Theta_{1}, \frac{\partial \Theta_{1}}{\partial t} \right)
\end{equation}

The power $q_u$ and the moment $m_u$ impacts are influenced by the high-speed pressure and therefore should depend on the velocity of the flow of $V$, the position of the feathers on the wing, their orientation and other factors associated with the adopted aerodynamic calculation scheme.

Let us assume that the feathers are completely rigid structural elements, and that a change in the feather's orientation does not influence the airflow around the remaining feathers. We also suppose that a wing as a whole keeps its laminar. Then, we obtain:
\begin{equation}\label{eq9}
q_{u} = \sum_{i}^{n(z)}{q_{u_{i}}}, ~~ m_{u} = \sum_{i}^{n(z)}{m_{u_{i}}},
\end{equation}
where $q_{u_{i}}$ and $m_{u_{i}}$ are the additional linear forces and the moment from the i-th feather, and the summation is carried out over all the feathers, covering the section frontier. The number of feathers is denoted as $n(z)$.

The tangent can be used to measure the rotation angle between the feathers on the upper surface and the wing profile $\beta_{i}\in[0,\beta^{-}],~\beta^{-}<0,$ and between the feathers on the lower surface and the wing profile $\beta_{i}\in[0,\beta^{+}],~\beta^{+}>0.$ 

\begin{figure}[thpb]
	\centering
	\includegraphics[scale=0.38]{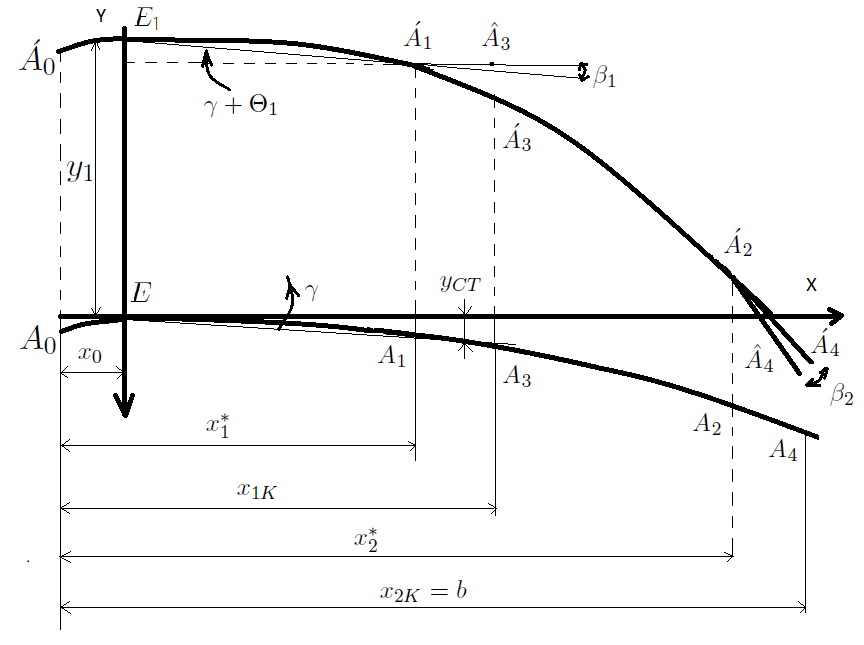}
	\caption{Wing profiles in a static position and during vibrations}
	\label{fig:4p}
\end{figure}

In Fig. \ref{fig:4p}:
\begin{itemize}
	\item $A_{0}~A_{1}~A_{2}~A_{4}$ are wing profiles (considered to be thin) in a static position (before vibrations);

	\item $\acute{A}_{0}~\acute{A}_{1}~\acute{A}_{2}~\acute{A}_{4}$ are wing profiles during oscillations; 

	\item $E$ and $E_{1}$ are SC's of the wing section in the static position and during vibrations, respectively;

	\item $X$-axis corresponds to the speed of the main stream;

	\item $Y$-axis is perpendicular to $X$-axis and to the axis of rigidity of the undeformed wing;

	\item $A_{1}~A_{3}~\acute{A}_{1}~\acute{A}_{3}$ are the front and back edges of the feather 1 in the neutral position on the corresponding profiles;

	\item $A_{2}~A_{4}~\acute{A}_{2}~\acute{A}_{4}$ are the front and back edges of the feather 2 (analog of the aileron) in the neutral position on relevant profiles;

	\item $x_{1}^{*}$ and $x_{1K}$ are the distances from the leading and trailing edges of the feather 1 to the leading edge of the wing;

	\item $x_{2}^{*}$ and $x_{2K}$ are similar parameters for the feather 2;

	\item $EE_{1}$ is the deflection of the wing;

	\item $\Theta_{1}$ is the angle of twisting of the wing near the point $E_{1}$;

	\item $\beta_{1}<0$ is the angle of deviation of the feather 1 from the neutral position;

	\item $\beta_{2}>0$ is the angle of deviation of the feather 2 from the neutral position;

	\item $\hat{A}_{3}$ and $\hat{A}_{4}$ are the trailing edges of the feathers 1 and 2, respectively, after their deviations.

\end{itemize}

Parameters $\Theta_{1}$, wing deflection $y_{CT}$ are ordinates in a static position, $\beta_{i}$ are considered small.

A similar figure for the aileron is presented in (\cite{A1} p. 143, Fig. 41).

Following the technique suggested by (\cite{A2}, pp.143-146), it can be shown that the influence of the $i$-th feather on the wing is generally calculated as follows:

\begin{equation}\label{eq10}
\begin{cases}
q_{u_{i}} = A_{i}V^{2}\beta_{i} + B_{i}V\dot{\beta}_{i}, \\
m_{u_{i}} = C_{i}V^{2}\beta_{i} + D_{i}V\dot{\beta}_{i},
\end{cases}
\end{equation}
where
$A_{i}=C_{y}^{\alpha}G_{i}\rho b^{2},~$
$B_{i}=C_{y}^{\alpha}H_{i}\rho b^{3},$

$C_{i}=-\left[ I_{i} + C_{y}^{\alpha}(\frac{x_{0}}{b}-\frac{1}{4})G_{i}\right]\rho b^{2},$

$D_{i}=-\left[ J_{i} + C_{y}^{\alpha}(\frac{x_{0}}{b}-\frac{1}{4})H_{i}\right]\rho b^{3},$

$G_{i}=\frac{1}{\pi}\left[ (\psi_{ik}-\psi_{i}^{*}) - (\sin{\psi_{ik}}-\sin{\psi_{i}^{*}})\right],$

\begin{multline}\notag
H_{i}=\frac{1}{2\pi}\left(
\cos{\psi_{i}^{*}}(\psi_{ik}-\psi_{i}^{*}) - (\sin{\psi_{ik}}-\sin{\psi_{i}^{*}})\right) 
-\cos{\psi_{i}^{*}}(\sin{\psi_{ik}}-\sin{\psi_{i}^{*}}) \\
+\frac{1}{2}\left(
(\psi_{ik}-\psi_{i}^{*}) +\frac{1}{2}(\sin{2\psi_{ik}}-\sin{2\psi_{i}^{*}})
\right);
\end{multline}

$I_{i}=\frac{1}{8}\left[ 2(\sin{\psi_{ik}}-\sin{\psi_{i}^{*}}) + (\sin{2\psi_{ik}}-\sin{2\psi_{i}^{*}})\right],$

\begin{multline}\notag
J_{i}=-\frac{1}{16}\left(
-2\cos{\psi_{i}^{*}}(\sin{\psi_{ik}}-\sin{\psi_{i}^{*}}) + (\psi_{ik}-\psi_{i}^{*})
\right)
-\frac{1}{16}\left(
(\frac{1}{2}-\cos{\psi_{i}^{*}})(\sin{2\psi_{ik}}-\sin{2\psi_{i}^{*}})
\right)
\\
-\frac{1}{16}\left( 
(\sin{\psi_{ik}}-\sin{\psi_{i}^{*}})
+\frac{1}{3}\left(
\sin{3\psi_{ik}}-\sin{3\psi_{i}^{*}}
\right)
\right);
\end{multline}

$x_{i}^{*}=\frac{b}{2}(1-\cos{\psi_{i}^{*}}),~~$
$x_{ik}=\frac{b}{2}(1-\cos{\psi_{ik}}),$

$\psi_{i}^{*}\in[0,\pi],~~$
$\psi_{ik}\in[0,\pi].$


According to \cite{A2,A3}, the solution of \eqref{sist6} near the flutter is given as
\begin{equation}\label{eq11}
\begin{cases}
y_{1}(z,t)=q(t)f(z),\\
\Theta_{1}(z,t)=r(t)\phi(z).
\end{cases}
\end{equation}

where $f(z)$ and $\phi(z)$ are vibration modes functions, satisfying the boundary conditions:
at $z=0,~f=0;~f'=0;~\phi=0;$
at $z=l,~f''=0;~f'''=0;~\phi'=0.$

Here for the sake of simplicity, we suggest $f'=\frac{\partial f}{\partial z}$,
$f'''=\frac{\partial^{3} f}{\partial z^{3}}.$

We substitute \eqref{eq11} into \eqref{sist6}, multiply the first equation by $f$ and the second equation by $\phi$ and then integrate from 0 to l. 
It gives after simple transformations with taking into account \eqref{eq9} and \eqref{eq10}
\begin{equation}\label{eq12}
\begin{cases}
a_{11}\ddot{q} + a_{12}\dot{q} + a_{13}q + b_{11}\ddot{r} + b_{12}\dot{r}+b_{13}r = Q(\beta,\dot{\beta}), \\
a_{21}\ddot{q} + a_{22}\dot{q} + b_{21}\ddot{r} + b_{22}\dot{r}+b_{23}r = M(\beta,\dot{\beta}),
\end{cases}
\end{equation}
where

$a_{11}=\int_{0}^{l}{mf^{2}dz},$

$a_{12}=C_{y}^{\alpha}\rho V\int_{0}^{l}{bfdz},$

$a_{13}=\int_{0}^{l}{\frac{d^{2}\left(EJf''\right)}{d{z}^{2}}fdz},$

$b_{11}=-\int_{0}^{l}{m\sigma_{T}f\phi dz},$

$b_{12}=-C_{y}^{\alpha}\rho V\int_{0}^{l}{\left(\frac{3}{4}b-x_{0}\right)bf\phi dz},$

$b_{13}=-C_{y}^{\alpha}\rho V^{2}\int_{0}^{l}{bf\phi dz},$

$a_{21}=\int_{0}^{l}{m\sigma_{T}f\phi dz}=-b_{11},$

$a_{22}=-C_{y}^{\alpha}\rho V\int_{0}^{l}{\left(x_{0}-\frac{b}{4}\right)bf\phi dz},$

$b_{21} = -\int_{0}^{l}{J_{m}\phi^{2}dz},$

$b_{22}=-\frac{\pi}{16}\rho V\int_{0}^{l}{b^{3}\phi^{2}dz} + 
C_{y}^{\alpha}\rho V\int_{0}^{l}{b(x_{0}-\frac{b}{4})\left(\frac{3}{4}b-x_{0}\right)\phi^{2} dz},$

$b_{23}=b_{23}^{(1)} + b_{23}^{(2)} = C_{y}^{\alpha}\rho V^{2}\int_{0}^{l}{b(x_{0}-\frac{b}{4})\phi^{2}dz} + 
\int_{0}^{l}{\frac{d\left(GJ_{k}\phi' \right)}{dz}\phi dz},$

$Q(\beta,\dot{\beta}) = \sum_{i=1}^{N}{\left(\bar{A}_{i}V^{2}\beta_{i} + \bar{B}_{i}V\dot{\beta}_{i}\right)},$

$\bar{A}_{i}=\int_{0}^{l}{A_{i}f dz},~~$
$\bar{B}_{i}=\int_{0}^{l}{B_{i}f dz},$

$M(\beta,\dot{\beta}) = \sum_{i=1}^{N}{\left( \bar{C}_{i}V^{2}\beta_{i} + \bar{D}_{i}V\dot{\beta}_{i}\right)},$

$\bar{C}_{i}=\int_{0}^{l}{C_{i}\phi dz},~~$
$\bar{D}_{i}=\int_{0}^{l}{D_{i}\phi dz},$

$\beta=col\left\{\beta_{i},~i={1,N} \right\},~~$
$\dot{\beta}=col\left\{\dot{\beta}_{i},~i={1,N} \right\},$
$N$ is the total feathers number.

Given the functions, $f$ and $\phi$ and the distributions of the mass and stiffness parameters of the wing (we consider time independent), 
the coefficients $a_{ij}$ and $b_{ij}, ~ i, j = 1,2,3$  can be figured out to be constants. 
The further results hardly rest on the choice of the functions $f$ and $\phi$. 
Note without going into details that these functions can be reasonably calculated, for example, by the successive approximations method. 
We complete \eqref{eq12} with the control equations

\begin{equation}\label{eq13}
\dot{\beta} = u,
\end{equation}
where
$u = col\left\{u_{i},~~i={1,N} \right\};~~$
$\beta_{i}\in[0,\beta^{+}], ~~ \beta^{+} > 0, ~~ i \in \bar{1,n^{+}},~$
,---, $n^{+}$ is the total number of feathers on the lower surface of the wing;

$\beta_{i}\in[\beta^{-},0], ~~ \beta^{-} < 0, ~~ i \in \bar{n^{+}+1,N},~$
,---, $n^{-}=N-n^{+}$ is total number of feathers on the upper surface of the wing.

Introduce:
\begin{equation}\label{eq14}
x=col\left\{q, \dot{q}, r, \dot{r} \right\} = col\left\{x_{i},~~i=\bar{1,4} \right\}.
\end{equation}

Then, we substitute \eqref{eq14} into \eqref{eq12} and reduce this system to the normal Cauchy form. After combining \eqref{eq12} and \eqref{eq13}, we get

\begin{equation}\label{eq15}
\begin{cases}
\dot{x}_{1}=x_{2},\\
\dot{x}_{2}=\sum_{k=1}^{4}{C_{1k}x_{k}} + F_{1}(\beta,u),\\
\dot{x}_{3}=x_{4},\\
\dot{x}_{4}=\sum_{k=1}^{4}{C_{2k}x_{k}} + F_{2}(\beta,u),\\
\dot{\beta}=u,
\end{cases}
\end{equation}
where
$F_{1}(\beta,u)=d_{11}Q+d_{12}M=\sum_{i=1}^{N}{\left(R_{1i}\beta_{i}+s_{1i}u_{i} \right)},$

$F_{2}(\beta,u)=d_{21}Q+d_{22}M=\sum_{i=1}^{N}{\left(R_{2i}\beta_{i}+s_{2i}u_{i} \right)},$

$R_{1i}=V^{2}\left( \bar{A}_{i}d_{11} + \bar{C}_{i}d_{12}\right),$

$s_{1i}=V\left( \bar{B}_{i}d_{11} + \bar{D}_{i}d_{12}\right),$

$d_{11}=\left[a_{11}\left(1-\frac{a_{21}b_{11}}{a_{11}b_{21}} \right) \right]^{-1},$

$d_{12}=-d_{11}{b_{11}} / {b_{12}},$

$R_{2i}=V^{2}\left( \bar{A}_{i}d_{21} + \bar{C}_{i}d_{22}\right),$

$s_{2i}=V\left( \bar{B}_{i}d_{21} + \bar{D}_{i}d_{22}\right),$

$d_{21}=-a_{21}\frac{d_{11}}{b_{21}},$

$d_{22}=\left(1 - a_{21}d_{12} \right) / b_{21},$

$C_{11}=-d_{11}a_{13},$

$C_{12}=-d_{11}(a_{12}-b_{11}\frac{a_{22}}{b_{21}}),$

$C_{13}=-d_{11}(b_{13}-b_{11}\frac{b_{23}}{b_{21}}),$

$C_{14}=-d_{11}(b_{12}-b_{11}\frac{b_{22}}{b_{21}}),$

$C_{21}=-a_{21}\frac{c_{11}}{b_{21}},$

$C_{22}=-\left( a_{22}+a_{21}c_{12}\right)/b_{21},$

$C_{23}=-\left( b_{23}+a_{21}c_{13}\right)/b_{21},$

$C_{24}=-\left( b_{22}+a_{21}c_{14}\right)/b_{21}.$

Now, we convert \eqref{eq7} using \eqref{eq11} and \eqref{eq14}
\begin{multline}\label{eq16}
E=\frac{1}{2}\int_{0}^{l}{mf^{2}dz}\dot{q}^{2} +
\frac{1}{2}\int_{0}^{l}{J_{m}\phi^{2}dz}\dot{r}^{2}
-\int_{0}^{l}{m\sigma_{T}f\phi dz} \dot{q}\dot{r} 
+\frac{1}{2}\int_{0}^{l}{EJ(f'')^{2}dz}q + \frac{1}{2}\int_{0}^{l}{GJ_{k}(\phi')^{2}dz}r = \\
\frac{1}{2}a_{13}x_{1}+\frac{1}{2}a_{11}x_{2}^{2}
-\frac{1}{2}b_{23}^{(2)}x_{3}-\frac{1}{2}b_{21}x_{4}^{2}-a_{21}x_{2}x_{4}\leq E_{*}.
\end{multline}

Finally, by integration by parts we get:
\begin{multline}\notag
\int_{0}^{l}{EJ(f'')^{2}dz} = \left. EJf''f' \right\rvert_{0}^{l} - 
\int_{0}^{l}{\left(EJf''\right)'f'dz} 
= 
-\frac{d\left( EJf''\right)}{dz}\left. f\right\rvert_{0}^{l} + \int_{0}^{l}{\frac{d^{2}\left( EJf''\right)}{dz^{2}}fdz} = a_{13},
\end{multline}

\begin{equation}\notag
\int_{0}^{l}{GJ_{k}\left(\phi' \right)^{2}dz} 
= GJ_{k}\phi'\left. \phi\right\rvert_{0}^{l} - \int_{0}^{l}{\frac{d\left(GJ_{k}\phi' \right)}{dz}\phi dz} 
= -b_{23}^{(2)}.
\end{equation}

Thus an equation \eqref{eq15} of the system is provided by a dynamics describing small flexural-twisting wing oscillations in a laminar flow taking into account the linear aerodynamic force and moment together with the equation \eqref{eq16} of limiting the total system energy. The changing rate of the inclination angle of the ``feather'' concerning the wing plane can be selected as a control parameter.

\section{CONTROL SYNTHESIS WITH THE SPEED-GRADIENT METHOD}

To equalize the forces impact on different parts of the wing, we use the Speed-Gradient Principle~\cite{A6, A7} to derive the feather angle control law.
According to this principle, all physical systems evolve along the shortest path in the direction of
thermodynamic equilibrium, which corresponds to the maximal value of entropy. In the Speed-Gradient algorithm, the maximal increment of entropy corresponds to the minimal value of the
energy in~\eqref{eq16}.

Up to this point, we have described the dynamic system taking into account how the solution should look in the optimal state, determined by \eqref{eq16}. However, a control law admitting the system to reach the desired state is absent. In this section such a law is produced.

First of all, we seek to control for the \eqref{eq15} with the criterion of \eqref{eq16}, using the speed-gradient method introduced in \cite{A6,A7}.

\begin{multline}\notag
\frac{dE}{dt} = \frac{1}{2}a_{13}x_{2}+a_{11}x_{2}\left( \sum_{k=1}^{4}{C_{1k}x_{k}} + F_{1}(\beta, u)\right)
-\frac{1}{2}b_{23}^{(2)}x_{4} 
- b_{21}x_{4}\left( \sum_{k=1}^{4}{C_{2k}x_{k}} + F_{2}(\beta, u)\right) \\
-a_{21}x_{4}\left( \sum_{k=1}^{4}{C_{1k}x_{k}} + F_{1}(\beta, u)\right) 
-a_{21}x_{2}\left( \sum_{k=1}^{4}{C_{2k}x_{k}} + F_{2}(\beta, u)\right),
\end{multline}

\begin{multline}\notag
\nabla_{u}\left( \frac{dE}{dt}\right) = 
col\left\{ (a_{11}x_{2}-a_{21}x_{4})\frac{\partial F_{1}}{\partial u_{i}}
-(b_{21}x_{4}+a_{21}x_{2})\frac{\partial F_{2}}{\partial u_{i}}, i={1,N}
\right\}=\\
col\left\{ (a_{11}x_{2}-a_{21}x_{4})s_{1i}
-(b_{21}x_{4}+a_{21}x_{2})s_{2i}, i={1,N}\right\}=\\
col\left\{ (a_{11}s_{1i}-a_{21}s_{2i})x_{2}
-(a_{21}s_{1i}+b_{21}s_{2i})x_{4}, i={1,N}\right\}=
col\left\{ \mu_{i}x_{2} + \nu_{i}x_{4}, i={1,N}\right\},
\end{multline}
where
$\mu_{i} = a_{11}s_{1i}-a_{21}s_{2i};~$
$\nu_{i} = -(a_{21}s_{1i}+b_{21}s_{2i}).$

Thus,
\begin{equation}\notag
\frac{du_{i}}{dt} = -\gamma_{i}(\mu_{i}x_{2}+\nu_{i}x_{4}), ~~ i={1,N}, ~~ \gamma_{i}>0
\end{equation}
\noindent or
\begin{equation}\notag
\frac{du_{i}}{dt} = -\gamma_{i}(\mu_{i}\dot{x}_{1}+\nu_{i}\dot{x}_{3}) \Rightarrow
u_{i} = -\gamma_{i}(\mu_{i}x_{1}+\nu_{i}x_{3}) + const_{i}.
\end{equation}

Since for $x_{1}=x_{2}=x_{3}=x_{4}=0$ all values of $u_{i}=0,$
and it means that $const_{i} = 0, i={1,N}.$

The resulting control equation derived in accordance with the Speed-Gradient Principle
\begin{equation}\label{eq17}
u_{i} = -\gamma_{i}(\mu_{i}x_{1}+\nu_{i}x_{3})
\end{equation}
is a control in the form of feedback on a deviation with constant coefficients.

\section{MULTI-AGENT CONTROL}

Utilization of the multiagent control allows for formation of emergent intelligence (intellectual resonance, swarm intelligence) or an occurrence of unexpected properties that a system possesses, but none of its individual elements has.

Each feather in the system aims to solve its own ``task'' of minimizing the deviation of the wing segment, to which the feather is attached, from its initial state. The action of each feather, in general, is not consolidated with other feathers actions, however a combination of all feathers impact results in a new property of the wing as the multiagent system of feathers to damp the vibrations. 

Now, consider the feathers as intelligent agents, so that each of them can receive information about the movement of the wing, exchange this information with other agents (transfer its information to them and receive their information), process the received data, and form local force and moment impacts on the wing, trying to keep the wing as close to its initial shape as possible (see Fig. \ref{fig:2p}). Assuming the feathers size is relatively small compared to the wing's surface, we relate a feather to some point on the wing's surface, when it is found in the neutral position.



Introduce for each feather ($i={1,N}$):
\begin{itemize}
	\item $\bar{z}_{i},\bar{\psi}_{i}$ --- coordinates of the point to which the $i$-th feather is related;
	\item $y_{1i}$ and $\Theta_{1i}$ --- deflection and angle of twisting of the wing at the location of the $i$-th feather (deviations from the curve \eqref{eq2});
	\item $N_{i}$ --- the set of feathers with which the $i$-th feather can exchange information;
	\item $\bar{b}_{ij}$ --- a non-negative weighting coefficient of the significance of information from the $i$-th feather to the $j$-th.
Here, we assume that $\bar{b}_{ij}=\bar{b}_{ji}$
and $\sum_{j\in N_{i}}{\bar{b}_{ij}}=1$ is the normalization condition;
	\item $\bar{b}_{ij}=0,$ if the i-th and $j$-th feathers are not connected informationally;
	\item $\bar{B}=[\bar{b}_{ij}]$ --- adjacency matrix.
\end{itemize}
 
According to \eqref{eq8} and \eqref{eq11}, for each $i={1,N}$
\begin{equation}\label{eq_Phi}
\begin{Vmatrix}  
w_{i}
\end{Vmatrix}
=
\begin{Vmatrix}  
y_{1i}\\ 
\dot{y}_{1i}\\
\Theta_{1i}\\
\dot{\Theta}_{1i}
\end{Vmatrix}
=
\begin{Vmatrix}  
qf(z_{i})\\ 
\dot{q}f(z_{i})\\
r\phi(z_{i})\\
\dot{r}\phi(z_{i})
\end{Vmatrix}
=
\begin{Vmatrix}  
\Phi_{i}\bar{x}
\end{Vmatrix}
< \epsilon
\end{equation}
for $t>t_{1}$ for an extended period of time;
$t_{1}$ --- the moment of reaching  $V_{flat}$;
$\Phi_{i}=diag\{f(z_{i}),f(z_{i}),\phi_{i}(z_{i}),\phi_{i}(z_{i})\}$

Moreover, we take into account that
\begin{equation}\notag
\left\|\Phi_{i}\bar{x} - \Phi_{j}\bar{x} \right\|
=\left\|\left[\Phi_{i} - \Phi_{j}\right]\bar{x} \right\|
\leq \left\|\Phi_{i}\bar{x}\right\| + \left\|\Phi_{j}\bar{x}\right\|
< 2\epsilon = \epsilon^{*}.
\end{equation}

The compensation for deviations from the stationary position is given in the model as
\begin{equation}\label{eq18}
L(\bar{x})=\frac{1}{2}\sum_{i=1}^{N}\sum_{j\in N_{i}}{\bar{b}_{ij}\left\|\left(\Phi_{i} - \Phi_{j}\right)\bar{x} \right\|^{2}}.
\end{equation}

The problem by analogy with \cite{A5} is formulated as follows.
In conditions of uniform rectilinear flight of the aircraft in a laminar flow, when approaching the critical speed of flexural-torsional flutter onset $V_{flat}$, it is required to find such controls $u_i$ for each feather in the system \eqref{eq15}, that would ensure the fulfillment of the target condition for the functional \eqref{eq18}:
\begin{equation*}
L(\bar{x}) \leq \epsilon^{*}
\end{equation*}

for a small given tolerance $\epsilon_{*}>0$ for $t>t_{1}$
during a sufficiently long period,
where $t_{1}$ is the time of reaching the critical flutter speed. The feather control laws are generated according to the Speed-Gradient Principle as discussed above.

\subsection{NON-MULTI-AGENT CONTROL SYNTHESIS}

Let us consider the functional
\begin{equation}\notag
L(x) = \frac{1}{2}\sum_{i=1}^{N}\sum_{j\in N_{i}}{\bar{b}_{ij}\left\|(\Phi_{i}-\Phi_{j})\bar{x}\right\|}^{2}\leq \epsilon^{*}.
\end{equation}

We suppose that $\bar{b}_{ij}\geq 0~\forall i,j$ which implies $L\geq 0.$

\begin{multline}\notag
\left\|(\Phi_{i}-\Phi_{j})\bar{x}\right\|^{2} = 
\left\|diag\{a,~b, ~c,~d \}\bar{x}\right\|^{2} = 
\left\|\left(f_{ij}x_{1},~f_{ij}x_{2},~\phi_{ij}x_{3},~\phi_{ij}x_{4} \right)^{T}\right\|^{2} = \\
f_{ij}^{2}(x_{1}^{2}+x_{2}^{2}) + \phi_{ij}^{2}(x_{3}^{2}+x_{4}^{2}),
\end{multline}
where 
$a = f(z_{i})-f(z_{j})$, $b = f(z_{i})-f(z_{j})$, $c = \phi(z_{i})-\phi(z_{j})$, $\Phi$ is defined in~\eqref{eq_Phi}, $d = \phi(z_{i})-\phi(z_{j})$ 
and 
$f_{ij}=f_{i}-f_{j}=f(z_{i})-f(z_{j}),$ $\phi_{ij}=\phi_{i}-\phi_{j}=\phi(z_{i})-\phi(z_{j}).$

So:
\begin{multline}\label{eq20}
L(x) = \frac{1}{2}\sum_{i=1}^{N}\sum_{j\in N_{i}}{\bar{b}_{ij}\left[f_{ij}^{2}(x_{1}^{2}+x_{2}^{2}) + \phi_{ij}^{2}(x_{3}^{2}+x_{3}^{2}) \right]} 
= \frac{1}{2}\left[\chi(x_{1}^{2}+x_{2}^{2}) + \lambda(x_{3}^{2}+x_{4}^{2}) \right],
\end{multline}
where $\chi = \sum_{i=1}^{N}\sum_{j\in N_{i}}{\bar{b}_{ij}f_{ij}^{2}}\geq 0$ and 
$\lambda = \sum_{i=1}^{N}\sum_{j\in N_{i}}{\bar{b}_{ij}\phi_{ij}^{2}}\geq 0$
are constants determined by the topology of the agent network (a wing in our model) (for given functions of the waveforms $f(z)$ and $\phi(z)$).
A singular situation appears once $\chi=0$ or $\lambda=0$.

Following the Speed-Gradient method of \cite{A6} from \eqref{eq20} and \eqref{eq15}, we obtain:
\begin{multline}\notag
\frac{dL}{dt}=\chi(x_{1}x_{2} + x_{2}\dot{x}_{2})+\lambda(x_{3}x_{4} + x_{4}\dot{x}_{4})= \\
\chi x_{2}\left[x_{1}+\sum_{k=1}^{4}{C_{1k}x_{k}} + F_{1}(\beta,u) \right] + 
\lambda x_{4}\left[x_{3}+\sum_{k=1}^{4}{C_{2k}x_{k}} + F_{2}(\beta,u) \right],
\end{multline}

\begin{equation}\notag
\nabla_{u}\dot{L} = col\left\{ \frac{\partial\dot{L}}{\partial u_{i}},~i=1,N\right\},
\end{equation}

\begin{multline*}
\frac{\partial}{\partial u_{i}}\left(\frac{dL}{dt}\right) = \chi x_{2}\frac{\partial F_{1}(\beta,u)}{\partial u_{i}}
+ \lambda x_{4}\frac{\partial F_{2}(\beta,u)}{\partial u_{i}} = \\
\chi x_{2}s_{1i} + \lambda x_{4}s_{2i} =
\chi s_{1i}\dot{x}_{1} + \lambda s_{2i}\dot{x}_{3}.
\end{multline*}

Consequently:
\begin{multline}\label{eq21}
\frac{d u_{i}}{dt} = -\gamma_{i}\left\{ \chi s_{1i}\dot{x}_{1} + \lambda s_{2i}\dot{x}_{3}\right\} \Rightarrow \\
\dot{\beta}_{i} = u_{i} = -\gamma_{i}\left\{ \chi s_{1i}{x}_{1} + \lambda s_{2i}{x}_{3}\right\}, ~ \gamma_{i}>0,~i=1,N,
\end{multline}
since the integration constant is zero for the same reasons as in~\eqref{eq17}.

The \eqref{eq17} and \eqref{eq21} have the same structure, but there are some differences. In fact, according to \eqref{eq15}, we can consider the coefficients $s_{1i}$ and $s_{2i}$ as the coefficients of influence of the $i$-th feather on the force factor in bending vibrations and on the moment factor in torsional vibrations, respectively. Actually, the values of these coefficients show the degree of participation of the $i$-th feather in wing dynamics.

In \eqref{eq17}, the feedback coefficients take into account the influence of the types of wing oscillations on each other, while in \eqref{eq21} the feedback coefficients for bending and torsional vibrations are strictly separated and determined through influence factors only on its type of oscillation. Equation \eqref{eq21} is not multiagent by nature since its dependence upon information about the state of other agents is static and it is invariant to the dynamics of the $i$-th feather.

\subsection{MULTI-AGENT CONTROL SYNTHESIS}
Now, we go back to \eqref{eq15}. We expand the vector of phase coordinates by introducing
\begin{equation}\label{eq22}
\tilde{x}_{i}=col\left\{ \Phi_{i}\bar{x},\beta_{i}\right\}.
\end{equation}

For this extended vector, we compose a functional analogously to \eqref{eq20}
\begin{equation}\label{eq23}
\tilde{L}=\frac{1}{2}\sum_{i=1}^{N}\sum_{j\in N_{i}}{\bar{b}_{ij}\left\| col\left\{ (\Phi_{i}-\Phi_{j})\bar{x}, \beta_{i}-\beta_{j}\right\}\right\|^{2}}.
\end{equation}
The application of the proposed approach can be justified by the fact that for minor deviations of the wing from the stationary position, determined by \eqref{eq2}, deviations of the feathers from their neutral position $\beta_{i},~i=1,N$ should be small as well. That is, at least for feathers on one side of the wing (lower/upper), we have that
\begin{equation}\notag
\left(\beta_{i}-\beta_{j}\right)^{2}\leq \beta_{i}^{2}+\beta_{j}^{2} < 2\epsilon^{2}_{\beta},
\end{equation}
where $\epsilon_{\beta}$ is a reasonably small number.

After simple transformations, we get
\begin{equation}\label{eq24}
\tilde{L}=L+\frac{1}{2}\sum_{i=1}^{N}\sum_{j\in N_{i}}{\bar{b}_{ij}\left(\beta_{i}-\beta_{j}\right)^{2}} < \epsilon^{*} +
\sum_{i=1}^{N}\sum_{j\in N_{i}}{\bar{b}_{ij}\epsilon^{2}_{\beta}} < \epsilon^{**}.
\end{equation}

Now, we apply the Speed-Gradient method:
\begin{equation}\notag
\frac{d\tilde{L}}{dt} = \dot{L} + \sum_{i=1}^{N}\sum_{j\in N_{i}}{\bar{b}_{ij}\left(\beta_{i}-\beta_{j}\right)\left(u_{i}-u_{j}\right)};
\end{equation}

\begin{equation}\notag
\nabla_{u}\left(\frac{d\tilde{L}}{dt}\right)=col\left\{ \frac{\partial\dot{L}}{\partial u_{i}} + 2\sum_{j\in N_{i}}{\bar{b}_{ij}(\beta_{i}-\beta_{j})}\right\}
\end{equation}

Finally, the control law is:
\begin{multline}\label{eq25}
\dot{\beta}_{i}=u_{i}=-\tilde{\gamma}_{i}\left(\chi s_{1i}\dot{x}_{1} + \lambda s_{2i}\dot{x}_{3}\right) 
-2\tilde{\gamma}_{i}\sum_{j\in N_{i}}{\bar{b}_{ij}(\beta_{i}-\beta_{j})},~ \tilde{\gamma}_{i}>0,~ i=1,N.
\end{multline}

It is important to pay attention to the multiagent nature of the control protocol of \eqref{eq22} and \eqref{eq25}, since the control signal for the rotation of each feather is formed on the basis of information about its own current state and the current state of the feathers connected with it. The connection is defined by the second term in \eqref{eq25} which states the dependence of feather $i$ angle adjustment $\dot{\beta}_i$ on the feather angle $\beta_i$ deviation from its neighbors' angles $\beta_j$, $j \in N_i$. At the same time, the first part of \eqref{eq25} describes control in the form of feedback with constant coefficients according to the speed of deviation of bending and torsional vibrations from the stationary state of (2). 

It is essential that if $\dot{x}_{1}=\dot{x}_{3}=0,$ then this does not entail $u_{i}=0,$ since in the general case there can be $x_{1}\neq 0$ and $x_{3}\neq 0.$ 
To reduce them, it is necessary to apply a control defined by \eqref{eq25}.
The expression $\dot{\beta}_{i}=u_{i}=0$ is true only in the case of complete absence of oscillations, when $x_{1} = x_{2} = x_{3} = x_{4} = 0,$ since only then $\beta_{i}=0,~p=1,N.$

The control of \eqref{eq25} does not explicitly depend on the time, at which the critical flutter speed is reached that allows using this control without any changes also to multiple transitions in speed across this boundary. Introduced control for each feather $i$ uses only local information about its own angle $\beta_i$ and angles of its nearest neighbors $\beta_j$. Suggested procedure, on the one hand, does not require data collection from all feathers on the wing to form the controls \eqref{eq25} and, on the other hand, requires only small data amount to compute the control inputs in \eqref{eq25} which are the essential features of multiagent approach.


\section{SIMULATION}
Let us compare the introduced control laws: synthesized with the speed-gradient method~\eqref{eq17}, non-multiagent~\eqref{eq22} and multiagent~\eqref{eq25} control laws.
In the simulation experiments the constants were chosen equal to the values (in some arbitrary units) listed below. Time step $\Delta_t = 10^{-5}$;
number of time instants $T = 10$; 
number of feathers $N = 5$;
air density $\rho = 1.225$; 
linear mass of the wing $m = 10$; 
wing section chord $b = 10$; 
wing length $l = 10$; 
from~\eqref{sist6} the derivative of the wing lift coefficient with respect to $\alpha$, $C^{\alpha}_y = 10$; 
airspeed $V = 10$; 
wing stiffness in bending $EJ = 50$; 
wing stiffness in torsion $GJ_k = 70$;
distance between stiffness centers and gravity centers in the wing cross section assumed to be constant and equal $\sigma_T = 0.1$; 
wing cross-section is assumed to be elliptical with height $a = 2$; 
linear mass moment of inertia of the wing relative to its stiffness axis $J_m = (\pi a b / 4) a^2 b^2 / (4 (a^2 + b^2))$; 
feather coordinates in Z axis $\bar{z} = l \cdot (0.1, 0.2,\ldots, 0.9)$;
feather coordinates in X axis according to formulas~\eqref{eq10} $\bar{\psi} = \pi / 4 \cdot (0.1, 0.2,\ldots, 0.9)$;
the coordinates of feathers trailing edges in X axis equal $x = 3/4 l$; 
the coordinates of feathers joint in X axis $x_\star = x - 0.01$;
the distance from the leading edge of the wing to the SC section $x_0 = l / 4$.

\begin{figure}
\includegraphics[scale=0.7]{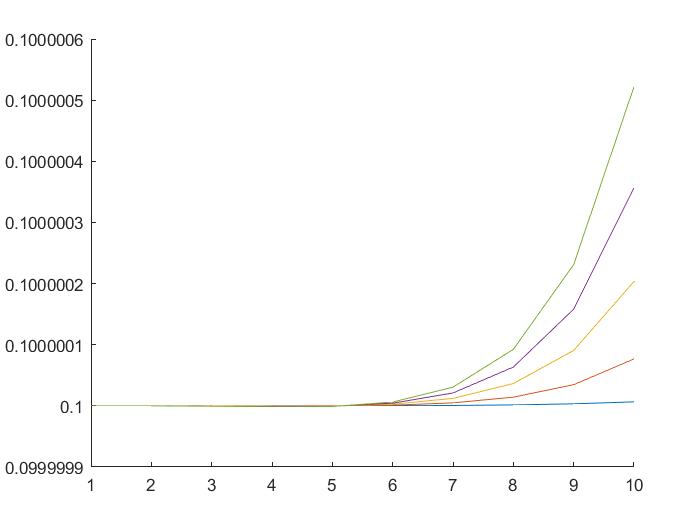}
\caption{feather angles $\beta$ change under the control synthesized with speed-gradient method}
\end{figure}

\begin{figure}
\includegraphics[scale=0.7]{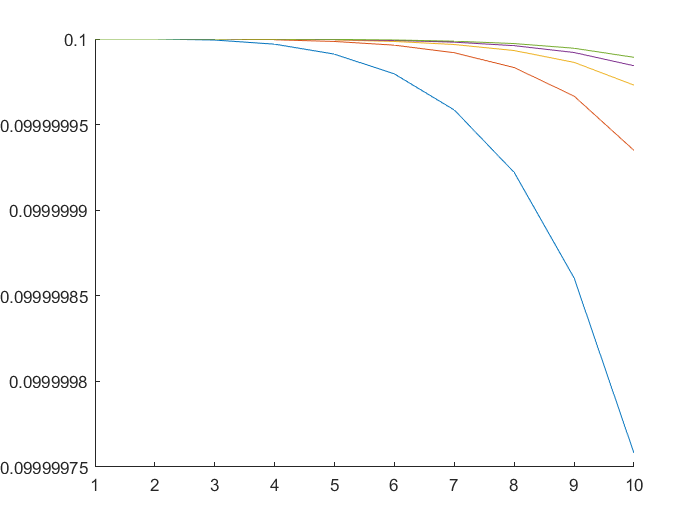}
\caption{feather angles $\beta$ change under the non-multiagent control law}
\label{fig:control2}
\end{figure}

\begin{figure}
\includegraphics[scale=0.7]{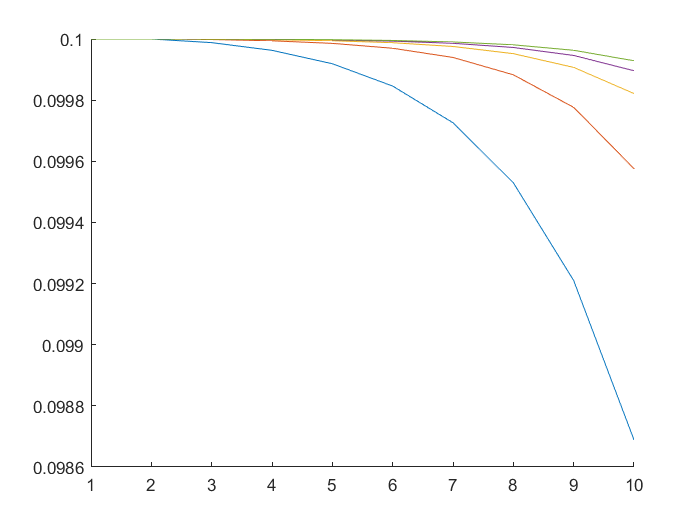}
\caption{feather angles $\beta$ change under the multiagent control law}
\label{fig:control3}
\end{figure}

The non-multiagent and multiagent control laws minimize the feather angles $\beta$ while the control law synthesized with the speed-gradient method tends to maximize the feather angles under given conditions. Multiagent and non-multiagent control treat the feathers separately and the difference between them lies in usage of the information from neighboring feathers in process of the control action synthesis. In case of multiagent control each feather ``takes into account'' the angles of its neighbors which allows to achieve faster control law. As could be seen from graphs in Figures~\ref{fig:control2} and~\ref{fig:control3} the multiagent control law (in Fig.~\ref{fig:control3}) minimizes the feather angles with a larger rate.

\section{CONCLUSIONS AND OUTLOOK}

This work is the first study of the authors, related to multiagent control of the wing with feathers aimed to avoid increasing wing oscillations when approaching flutter.
In the article
- a mathematical model of the bending-torsional vibrations of an airplane wing with controlled feathers on its surface is given;
- three different statements of the control problem are considered, which differ by the goal functional; 
- the three control laws \eqref{eq17}, \eqref{eq21} and \eqref{eq25} are synthesized by the Speed-Gradient method. Only one of them: \eqref{eq25} is multiagent.

\eqref{eq21} is an "intermediate" one for the synthesis of a multiagent control law. It has a similar structure to \eqref{eq17}, but takes into account the presence of other feathers and their contribution.
However, the information about them remains static. It means that the state and dynamics of other pens is not considered. The multiagent control law allows for each feather to take into account information about its own current state and about the current state of feathers in the area of the wing where it is located. As a rule, this allows for more precise tuning and quicker tuning to external factors, which makes this equation the most promising.

In the future, we advise to study the effectiveness of the obtained control laws and to compare them. The most critical indicator in the comparison should be the time, during which the system is able to damp vibrations to a safe level and hold them. The relevance of this indicator is due to the rather fast process of increasing wing oscillations during flutter.
Another promising area for the further research is the development of multiagent control of feathers following the example of a swarm and research its effectiveness.

\section*{ACKNOWLEDGMENT}

Sections 1-5 of the work were supported by the IPME RAS by Russian Science Foundation (project no. 21-19-00516).
Section 6 was supported by the St. Petersburg State University (project No. 73555239).

\end{document}